\documentclass[preprint,1p,times]{elsarticle}

\usepackage{amssymb}
\usepackage{amsthm}
\usepackage{amsmath}

\usepackage{lineno}

\journal{Linear Algebra and Its Applications}
\newtheorem{thm}{Theorem}[section]
\newtheorem{cor}[thm]{Corollary}
\newtheorem{lem}[thm]{Lemma}

\newtheorem{prop}[thm]{Proposition}

\newtheorem{rem}[thm]{Remark}

\newproof{pf}{Proof}

\begin{document}
\begin{frontmatter}



\title{Characterization of the oblique projector $U(VU)^{\dag}V$ with application to constrained least squares\tnoteref{label1}}%
\tnotetext[label1]{I would like to thank Adi Ben-Israel, Laura Rebollo-Neira, G\"{o}tz Trenkler
and two anonymous referees for helpful comments and pointers to references.
Any errors in the manuscript are my responsibility.}

\author{Ale\v{s} \v{C}ern\'{y}}

\address{Cass Business School, City University London}

\begin{abstract}
We provide a full characterization of the oblique projector $U(VU)^{\dag }V$
in the general case where the range of $U$ and the null space of $V$ are not
complementary subspaces. We discuss the new result in the context of
constrained least squares minimization which finds many applications in
engineering and statistics.\\
\vspace*{-6pt}\vspace*{-13pt}\\
\end{abstract}


\end{frontmatter}



\section{Introduction}

Let $E\in \mathbb{C}^{m\times m}$ be idempotent, $E^{2}=E$. The null space
and range of any idempotent matrix are complementary, cf. \cite[Theorem 2.8]%
{ben-israel.greville.03}, 
\[
R(E)+N(E)=\mathbb{C}^{m},R(E)\cap N(E)=\{0\},
\]%
and we say that $E$ is an oblique projector onto $R(E)$ along $N(E)$. For
any two complementary subspaces of $\mathbb{C}^{m}$ we denote the oblique
projector onto $L$ along $M$ by $P_{L,M}$. The orthogonal projector onto $L$
is denoted by $P_{L}:=P_{L,L^{\bot }}$, where $L^{\bot }$ is the orthogonal
complement of $L.$ Oblique projectors arise in numerous engineering and
statistical applications, see \cite[Chapter 8]{ben-israel.greville.03}, \cite%
{corach.maestripieri.05} and references therein. Many of their properties
follow from the general solution to the matrix equation $XAX=X$ studied in
1960-ies in the context of the various pseudoinverses, cf. \cite{trenkler.94}%
. This literature is mature, with excellent monographs such as \cite%
{ben-israel.greville.03}. In particular it is very well understood how to
construct an oblique projector with a prescribed range and null space.

\begin{prop}
\label{prop: intro}Let $L,M$ be complementary subspaces of $\mathbb{C}^{m}$.
For any two matrices $U,V$ with $R(U)=L$ and $N(V)=M$ one has%
\[
P_{L,M}=U(VU)^{\dag }V,
\]%
where the superscript \textquotedblleft \thinspace $\dag $%
\/\textquotedblright\ denotes the Moore-Penrose inverse. If $U$ and $V$ are
in addition orthogonal projectors (i.e. they are Hermitian and idempotent)
one obtains an even simpler form due to Greville \cite[(3.1) and Theorem 2]%
{greville.74}, 
\begin{equation}
P_{L,M}=P_{L}(P_{M^{\bot }}P_{L})^{\dag }P_{M^{\bot }}=(P_{M^{\bot
}}P_{L})^{\dag }.  \label{eq: Greville74}
\end{equation}
\end{prop}

The converse problem of characterizing the range and null space of a given
idempotent matrix has not received the same amount of attention. The
motivation for studying idempotents of the form $U(VU)^{\dag }V$ in the
general case where $R(U)+N(V)\subsetneq \mathbb{C}^{m}$ and/or $R(U)\cap
N(V)\neq \{0\}$ comes, among others, from constrained least squares optimization with a
range of applications mentioned above. Briefly, the problem%
\[
\min_{x\in \mathbb{C}^{n}}\left\Vert A_{1}x-b_{1}\right\Vert ^{2},\textrm{
subject to }A_{2}x=b_{2},
\]
gives rise to the projector $D_{2}(A_{1}D_{2})^{\dag }A_{1}$ where $D_{2}\ $%
is an arbitrary but fixed matrix with the property $R(D_{2})=N(A_{2}).$ In
this situation we typically have neither $R(D_{2})+N(A_{1})=\mathbb{C}^{m}$
nor $R(D_{2})\cap N(A_{1})=\{0\}$. Oblique projectors of the
form $U(VU)^{\dag }V$ with $R(U)+N(V)=\mathbb{C}^{m}$ and $R(U)\cap N(V)\neq
\{0\}$ feature also in signal reconstruction, cf. \cite{hirabayashi.unser.07}%
.

Given that $U(VU)^{\dag }V$ has a wide range of applications it is desirable to
understand its geometric nature. One might conjecture that in general%
\begin{align}
&U(VU)^{\dag }V=P_{L,M},\textrm{ where}\label{eq: conj1} \\ 
&L=P_{R(U)}N(V)^{\bot }=R(U)\cap (R(U)\cap N(V))^{\bot }, \\ 
&M=N(V)+(N(V)+R(U))^{\bot },\label{eq: conj3}
\end{align}%
but the behaviour of the projector is somewhat more intricate and cannot be
described based on the knowledge of $R(U)$ and $N(V)$ alone. The conjecture (%
\ref{eq: conj1})-(\ref{eq: conj3}) turns out to be true only when both $U$ and $V$ are
orthogonal projectors. Surprisingly, the main tool in proving the general
result is the Zlobec formula \cite{zlobec.70} in conjunction with
Proposition \ref{prop: intro}.

The result presented here is different from the problem discussed by Rao and
Yanai \cite{rao.yanai.79} in which projectors onto and along two given
subspaces are considered under the assumption that the subspaces are not
necessarily spanning the whole space. In such a situation, the projector no
longer needs to be idempotent.

The paper is organized as follows. In section 2 we introduce required
terminology and notation, we establish the main tools and prove Proposition %
\ref{prop: intro}. In section 3 we state and prove the main result. In
section 4 we discuss application of the main result to constrained least
squares minimization and the link to the minimal norm solution of Eld\'{e}n 
\cite{elden.80}.

\section{Preliminaries}

We use notation of \cite{ben-israel.greville.03}. $A^{\ast }$ denotes the
conjugate transpose of matrix $A.$ We write $r(A),$ $R(A),$ $N(A)$ for the rank,
range and null space of $A,$ respectively. Consider the following relations%
\begin{align}
AXA &=A,  \tag*{(I.1)} \\
XAX &=X,  \tag*{(I.2)} \\
AX &=(AX)^{\ast },  \tag*{(I.3)} \\
XA &=(XA)^{\ast }. \tag*{(I.4)}
\end{align}%
We write $X\in A\{i,j,\,\ldots ,k\}$, if $X$ satisfies conditions (I.$i$), (I.$j$), $\ldots,$
(I.$k$). $A^{\dag }$ denotes the Moore-Penrose inverse which is the unique
element of $A\{1,2,3,4\}.$ The following theorem is our main tool.

\begin{thm}[{\protect\cite[Theorem 2.13]{ben-israel.greville.03}}]
\label{th: main}Let $A\in \mathbb{C}^{m\times n},\tilde{U}\in \mathbb{C}%
^{n\times s},\tilde{V}\in \mathbb{C}^{t\times m}$ and 
\[
Z=\tilde{U}(\tilde{V}A\tilde{U})^{(1)}\tilde{V}, 
\]%
where $(\tilde{V}A\tilde{U})^{(1)}$ is a fixed but arbitrary element of $(%
\tilde{V}A\tilde{U})\{1\}$. Then

a) $Z\in A\{1\}$ if and only if $r(\tilde{V}A\tilde{U})=r(A);$

b) $Z\in A\{2\}$ and $R(Z)=R(\tilde{U})$ if and only if $r(\tilde{V}A\tilde{U%
})=r(\tilde{U});$

c) $Z\in A\{2\}$ and $N(Z)=N(\tilde{V})$ if and only if $r(\tilde{V}A\tilde{U%
})=r(\tilde{V});$

d) $Z=A_{R(\tilde{U}),N(\tilde{V})}^{(1,2)}$ if and only if $r(\tilde{U})=r(%
\tilde{V})=r(\tilde{V}A\tilde{U})=r(A),$ where $A_{R(\tilde{U}),N(\tilde{V}%
)}^{(1,2)}$ is the unique element of $A\{1,2\}$ with range $R(\tilde{U})$
and null space $N(\tilde{V})$, also known as the oblique pseudoinverse (cf. 
\cite{milne.68}).
\end{thm}

\begin{cor}
The Zlobec formula \cite{zlobec.70}, 
\begin{equation}
A^{\dag }=A^{\ast }(A^{\ast }AA^{\ast })^{(1)}A^{\ast },  \label{eq: Zlobec}
\end{equation}
is now obtained by setting $\tilde{U}=\tilde{V}=A^{\ast }$ in part d) and
arguing $A_{R(A^{\ast }),N(A^{\ast })}^{(1,2)}=A^{\dag }$.
\end{cor}

\noindent The following is a pre-cursor to the main result in this note. The
\textquotedblleft if\textquotedblright\ part appears, for example, in \cite[%
(3.51)]{pollock.79}.

\begin{cor}
\label{cor: main}$\tilde{U}(\tilde{V}\tilde{U})^{(1)}\tilde{V}=P_{R(\tilde{U}%
),N(\tilde{V})}$ if and only if $r(\tilde{V}\tilde{U})=r(\tilde{V})=r(\tilde{%
U}).$
\end{cor}

Next we show that the form $U(VU)^{\dag }V$ covers all idempotent matrices.

\begin{lem}
\label{lem: complementarity}Let $U\in \mathbb{C}^{m\times p},V\in \mathbb{C}%
^{q\times m}.$ $R(U)$ and $N(V)$ are complementary subspaces of $\mathbb{C}%
^{m}$ if and only if $r(U)=r(V)=r(VU)$.
\end{lem}

\begin{pf}
If: By Corollary \ref{cor: main} $U(VU)^{\dag }V=P_{R(U),N(V)}$ which
implies that $R(U),N(V)$ are complementary.

Only if: i) complementarity implies $\dim (R(U))+\dim (N(V))=m.$ On
rearranging we obtain $r(U)=m-\dim (N(V))$ and by the rank-nullity theorem $%
r(U)=r(V).$

ii) Complementarity also implies $R(U)\cap N(V)=\{0\}$ which yields $%
N(VU)=N(U)$. By rank-nullity theorem we obtain $r(VU)=r\left( U\right)
.\hfill\square $
\end{pf}

\begin{prop}
\label{prop: Eform}Matrix $E\in \mathbb{C}^{m\times m}$ is idempotent if and
only if there are matrices $U\in \mathbb{C}^{m\times p},V\in \mathbb{C}%
^{q\times m}$ such that 
\begin{equation}
E:=U(VU)^{\dag }V.  \label{eq: Eform}
\end{equation}
\end{prop}

\begin{pf}
The `if' statement follows easily from (\ref{eq: Eform}) and (I.2), $$E^{2}=U(VU)^{\dag
}VU(VU)^{\dag }V=E.$$ The `only if' part: construct $U$ so that its columns
form a basis of $R(E)$ and construct $V^{\ast }$ so that its columns form
the basis of $N(E)^{\bot }.$ This implies $R(U)=R(E),N(V)=N(E)$. Since $E$
is idempotent $R(U),N(V)$ are by construction complementary and from Lemma %
\ref{lem: complementarity} we obtain $r(U)=r(V)=r(VU).$ By Corollary \ref%
{cor: main} $U(VU)^{\dag }V=P_{R(E),N(E)}=E.\hfill\square $
\end{pf}

\begin{rem}
A comprehensive characterization of projectors appears in \cite{trenkler.94}%
. Proposition \ref{prop: Eform} resembles a result of Mitra \cite[Theorem 3a]%
{mitra.68} who shows that all idempotent matrices are of the form $\tilde{U}(%
\tilde{V}\tilde{U})^{(1,2)}\tilde{V}$ where $(\tilde{V}\tilde{U})^{(1,2)}$
is an arbitrary element of $\tilde{V}\tilde{U}\{1,2\}$. This result is
generalized further in \cite[Theorem 2.13]{ben-israel.greville.03} to the
form $\tilde{U}(\tilde{V}\tilde{U})^{(1)}\tilde{V}$, see Corollary \ref{cor:
main}. Proposition \ref{prop: Eform} goes in the opposite direction in order
to avoid the ambiguity associated with $\{1,2\}$-inverses.
\end{rem}

To conclude we provide a proof of Proposition \ref{prop: intro}.

\begin{pf}[Proposition \protect\ref{prop: intro}]
The first statement follows from the `only if' part in the proof of
Proposition \ref{prop: Eform}. The second part follows from identities $%
(P_{M^{\bot }}P_{L})^{\dag }=P_{L}(P_{M^{\bot }}P_{L})^{\dag }=(P_{M^{\bot
}}P_{L})^{\dag }P_{M^{\bot }},$ see \cite[Exercise 2.57]%
{ben-israel.greville.03}.\hfill$\square $
\end{pf}

\section{Result}

\begin{thm}
\label{th: result}Given two arbitrary matrices $U\in \mathbb{C}^{m\times
p},V\in \mathbb{C}^{q\times m}$ the matrix $E=U(VU)^{\dag }V$ is idempotent
with range and null space given by%
\begin{align}
R(E) &= R(UU^{\ast }V^{\ast })=R(UU^{\ast }V^{\ast }V) = R(U)\cap ((UU^{\ast })^{\dag }(R(U)\cap N(V)))^{\bot },%
\label{eq: result1} \\
N(E) &=  N(U^{\ast }V^{\ast }V)=N(UU^{\ast }V^{\ast }V) = N(V)\oplus (V^{\ast }V)^{\dag }(R(U)+N(V))^{\bot }.
\label{eq: result2}%
\end{align}
\end{thm}

\begin{pf}
By Zlobec's formula (\ref{eq: Zlobec}) with $A=$ $VU$ we obtain 
\[
E=UU^{\ast }V^{\ast }(U^{\ast }V^{\ast }VUU^{\ast }V^{\ast })^{(1)}U^{\ast
}V^{\ast }V.
\]%
Setting $\tilde{U}=UU^{\ast }V^{\ast },\tilde{V}=U^{\ast }V^{\ast }V$ we
claim $r(\tilde{U})=r(\tilde{V})=r(\tilde{V}\tilde{U})=r(VU).$ Indeed,
\begin{align}
r(VU) & =  r(VUU^{\ast }V^{\ast })=r(VUU^{\ast }V^{\ast }VUU^{\ast }V^{\ast
}) \leq  r(U^{\ast }V^{\ast }VUU^{\ast }V^{\ast })=r(\tilde{V}\tilde{U}), \label{eq: inequalities1}\\ 
r(\tilde{V}\tilde{U}) & \leq  r(\tilde{U})=r(UU^{\ast }V^{\ast })\leq r(U^{\ast }V^{\ast })=r(VU), \\ 
r(\tilde{V}\tilde{U}) & \leq  r(\tilde{V})=r(U^{\ast }V^{\ast }V)\leq
r(U^{\ast }V^{\ast })=r(VU).\label{eq: inequalities3}
\end{align}%
Corollary \ref{cor: main} yields $R(E)=R(\tilde{U}),N(E)=N(\tilde{V}).$ From%
\[
r(VU) =r(VUU^{\ast }V^{\ast })=r(VUU^{\ast }V^{\ast }VUU^{\ast }V^{\ast })
\leq r(UU^{\ast }V^{\ast }V)\leq r(U^{\ast }V^{\ast })=r(VU),
\]
and from (\ref{eq: inequalities1})-(\ref{eq: inequalities3}) we obtain $r(VU)=r(UU^{\ast }V^{\ast
})=r(UU^{\ast }V^{\ast }V)$ which implies $R(UU^{\ast }V^{\ast })=R(UU^{\ast
}V^{\ast }V)$. The proof of $N(U^{\ast }V^{\ast }V)=N(UU^{\ast }V^{\ast }V)$
proceeds similarly by showing $r(U^{\ast }V^{\ast }V)=r(UU^{\ast }V^{\ast
}V) $.

To show the last equality in (\ref{eq: result2}) we observe $\mathbb{C}%
^{m}=N(V)\oplus R(V^{\ast })$. Since $N(V)\subseteq N(U^{\ast }V^{\ast }V)$
we have 
\begin{equation}
N(U^{\ast }V^{\ast }V)=N(V)\oplus (R(V^{\ast })\cap N(U^{\ast }V^{\ast }V)).
\label{eq: N(UVV)}
\end{equation}%
Continuing with the second term on the right hand side we obtain%
\begin{align*}
y\in R(V^{\ast })\cap N(U^{\ast }V^{\ast }V) &\iff (V^{\ast }Vy\in
N(U^{\ast })\cap R(V^{\ast }))\wedge (y\in R(V^{\ast })) \\
&\iff y\in (V^{\ast }V)^{\dag }(N(U^{\ast })\cap R(V^{\ast })),
\end{align*}%
which yields%
\begin{equation}
R(V^{\ast })\cap N(U^{\ast }V^{\ast }V)  =  (V^{\ast }V)^{\dag}(R(U)^{\bot }\cap N(V)^{\bot })
 =  (V^{\ast }V)^{\dag }(R(U)+N(V))^{\bot }. \label{eq: N(UVV)2}
\end{equation}%
On substituting (\ref{eq: N(UVV)2}) into (\ref{eq: N(UVV)}) we obtain the
desired result.

The last equality in (\ref{eq: result1}) is obtained by writing $R(UU^{\ast
}V^{\ast })=N(VUU^{\ast })^{\bot }$ and then evaluating $N(VUU^{\ast })$ by
exchanging the role of $U$ and $V^{\ast }$ in (\ref{eq: N(UVV)}) and (\ref%
{eq: N(UVV)2}).\hspace*{\fill} $\square $
\end{pf}

\begin{rem}
\label{rem: result}Special cases of Theorem \ref{th: result} include
situations covered by Corollary \ref{cor: main} in which $r(U)=r(V)=r(VU)$
and we have $R(E)=R(U),N(E)=N(V);$ the Langenhop form \cite[Lemma 2.2]%
{langenhop.67} with $VU=I$ is a case in point. The Greville formula (\ref%
{eq: Greville74}) also falls into this category. Hirabayashi and Unser \cite[%
Lemma 3]{hirabayashi.unser.07} encounter the case $R(U)+N(V)=\mathbb{C}^{m}$
and $R(U)\cap N(V)\neq \{0\},$ yielding $R(E)=R(UU^{\ast }V^{\ast
}),N(E)=N(V)$.
\end{rem}

\section{Application}

\begin{prop}
\label{prop: ls}Let $A_{1}\in \mathbb{C}^{m\times n},b_{1}\in \mathbb{C}%
^{m},A_{2}\in \mathbb{C}^{k\times n},r(A_{2})=k\geq 1,b_{2}\in \mathbb{C}%
^{k} $. Solutions of the problem 
\begin{equation}
\min_{x\in \mathbb{C}^{n}}\left\Vert A_{1}x-b_{1}\right\Vert ^{2},\textrm{
subject to }A_{2}x=b_{2},  \label{eq: cls}
\end{equation}%
lie in the set%
\begin{equation}
\Xi =\{D_{2}(A_{1}D_{2})^{\dag }A_{1}A_{1}^{\dag
}b_{1}+(I-D_{2}(A_{1}D_{2})^{\dag }A_{1})(A_{2}^{\dag }b_{2}+z):z\in
N(A_{2})\},  \label{eq: clsformula}
\end{equation}%
where $D_{2}\ $is an arbitrary but fixed matrix with the property $%
R(D_{2})=N(A_{2}).$
\end{prop}

\begin{pf}
See \cite[Exercise 3.10]{ben-israel.greville.03}.\hfill $\square$
\end{pf}

In general, the projector $D_{2}(A_{1}D_{2})^{\dag }A_{1}$ will depend on how $D_{2}$ is chosen.
However, Theorem \ref{th: result} shows that there is a special case when 
$D_{2}(A_{1}D_{2})^{\dag }A_{1}$ is actually invariant to the choice of $D_2$. 

\begin{cor}
\label{cor : new}Using the notation of Proposition \ref{prop: ls} assume
further $r(A_{1})=n$. Then%
\[
D_{2}(A_{1}D_{2})^{\dag }A_{1}=P_{N(A_{2}),(A_{1}^{\ast
}A_{1})^{-1}R(A_{2}^{\ast })},
\]%
and $\Xi $ is a singleton,
\[
\Xi=\{A_{1}^{\dag }b_{1}+(A_{1}^{\ast }A_{1})^{-1}A_{2}^{\ast }(A_2 (A_{1}^{\ast }A_{1})^{-1}A_{2}^{\ast })^{-1}
(b_{2}-A_{2}A_{1}^{\dag}b_{1})\}.
\]
\end{cor}

\begin{pf}
We have $N(A_{1})=0$ and by Theorem \ref{th: result} 
\begin{align*}
R(D_{2}(A_{1}D_{2})^{\dag }A_{1}) &=R(D_{2})\cap (\{0\})^{\bot }=N(A_{2}),
\\
N(D_{2}(A_{1}D_{2})^{\dag }A_{1}) &=(A_{1}^{\ast }A_{1})^{-1}R(D_{2})^{\bot
}=(A_{1}^{\ast }A_{1})^{-1}R(A_{2}^{\ast }).
\end{align*}%
This implies $(I-D_{2}(A_{1}D_{2})^{\dag }A_{1})z=0$ for all $z\in N(A_{2})$ and by Proposition \ref{prop: intro}
$$(I-D_{2}(A_{1}D_{2})A_1)= (A_{1}^{\ast }A_{1})^{-1}A_{2}^{\ast }(A_2 (A_{1}^{\ast }A_{1})^{-1}A_{2}^{\ast })^{-1}A_2.$$ 
The rest follows from Proposition \ref{prop: ls}.\hfill $\square$
\end{pf}

Note that Corollary \ref{cor : new} is not covered by Corollary \ref{cor:
main} since $n-k=r(D_{2})=r(A_{1}D_{2})<r(A_{1})=n$. In situations where the
choice of $D_{2}$ impacts on the projector $D_{2}(A_{1}D_{2})^{\dag }A_{1}$
Theorem \ref{th: result} guides us to the convenient choice of $D_{2}$ which
simplifies the geometry of the result and also helps to identify the element 
of $\Xi$ with minimal distance from a given reference point.

\begin{cor}
\label{cor: new2}Using the notation of Proposition \ref{prop: ls} the following statements hold:
\begin{enumerate}
\item The constrained least squares minimizer in (\ref{eq: cls}) lies in the set 
\begin{equation}
\Xi =\{A_{1}^{\dag }b_{1}+P_{\mathcal{Y},\mathcal{%
X}}(A_{2}^{\dag }b_{2}-A_{1}^{\dag }b_{1})+z:z\in N(A_{1})\cap N(A_{2})\},
\label{eq: Xi_result}
\end{equation}
with%
\begin{align}
P_{\mathcal{X},\mathcal{Y}}&=I-P_{\mathcal{Y},\mathcal{X}}=(A_1 (I-A_{2}^{\dag }A_{2}))^\dag A_1,\label{eq: PYX}\\
\mathcal{X} &=P_{N(A_{2})}R(A_{1}^{\ast })=N(A_{2})\cap
(N(A_{2})\cap N(A_{1}))^{\bot },  \label{eq: calX} \\
\mathcal{Y} &=N(A_{1})\oplus(A_{1}^{\ast }A_{1})^{\dag }(N(A_{1})+
N(A_{2}))^\bot.  \label{eq: calY}
\end{align}
\item The element of $\Xi$ with the smallest Euclidean norm is given by 
\[
\xi:=A_{1}^{\dag }b_{1}+P_{\mathcal{Y},\mathcal{X}}(A_{2}^{\dag }b_{2}-A_{1}^{\dag }b_{1}).
\]
\item For any $y\in\mathbb{C}^n$ the solution of $\min_{x\in\Xi}||x-y\,||$ is given by 
\begin{equation}
\psi(y):=\xi+P_{N(A_1)\cap N(A_2)}y.
\label{eq: psidef}
\end{equation}
\end{enumerate}
\end{cor}

\begin{pf}
1. On setting $D_{2}=I-A_{2}^{\dag }A_{2}=P_{N(A_2)}$ Proposition \ref{prop: ls} and Theorem \ref{th:
result} yield%
\begin{equation}
\Xi = A_{1}^{\dag }b_{1}+P_{\mathcal{Y},\mathcal{X}}(A_{2}^{\dag }b_{2}-A_{1}^{\dag }b_{1}+N(A_2)),
\label{eq: Xi_aux}
\end{equation}
with $P_{\mathcal{X},\mathcal{Y}}$, $\mathcal{X}$ and $\mathcal{Y}$ given in (\ref{eq: PYX})-(\ref{eq: calY}). From (\ref{eq: calX}) we obtain 
$N(A_2)=\mathcal{X}\oplus (N(A_1)\cap N(A_2))$ which implies 
\begin{equation}
P_{\mathcal{Y},\mathcal{X}}N(A_2)=P_{\mathcal{Y},\mathcal{X}}(N(A_1)\cap N(A_2))=N(A_1)\cap N(A_2),
\label{eq: PYX(NA2)}
\end{equation}
the last equality following from $N(A_1)\cap N(A_2)\subseteq\mathcal{Y}$. Substitution of (\ref{eq: PYX(NA2)}) into (\ref{eq: Xi_aux}) 
yields (\ref{eq: Xi_result}).
\medskip\\
\noindent 2. By (\ref{eq: calX}) we have $\mathcal{X}\subseteq (N(A_{2})\cap N(A_{1}))^{\bot }=R(A_1^\ast)+R(A_2^\ast)$. Consequently
\begin{equation}
P_{\mathcal{Y},\mathcal{X}}(R(A_1^\ast)+R(A_2^\ast))=(I-P_{\mathcal{X},\mathcal{Y}})(R(A_1^\ast)+R(A_2^\ast))\subseteq R(A_1^\ast)+R(A_2^\ast).
\label{eq: PXYrange}
\end{equation} 
This implies 
\begin{equation}
\xi\in R(A_1^\ast)+R(A_2^\ast)=(N(A_{2})\cap N(A_{1}))^{\bot }. 
\label{eq: xi_orth}
\end{equation}
By (\ref{eq: Xi_result}) $x-\xi\in N(A_1)\cap N(A_2)$ for any $x\in\Xi$ which together with (\ref{eq: xi_orth}) yields 
\[
\|x\,\|^2=\|x-\xi+\xi\,\|^2=\|x-\xi\,\|^2+\|\xi\,\|^2\textrm{ for all } x\in\Xi.
\]
\noindent 3. By (\ref{eq: Xi_result}), (\ref{eq: psidef}) and (\ref{eq: xi_orth}) we obtain $x-\psi(y)\in N(A_{2})\cap N(A_{1})$ and $\psi(y)-y\in (N(A_{2})\cap N(A_{1}))^{\bot }$ which implies
$
\|x-y\,\|^2=\|x-\psi(y)+\psi(y)-y\,\|^2=\|x-\xi\,\|^2+\|\xi-y\,\|^2,
$
for all $x\in\Xi$.\hfill$\square$
\end{pf}

\begin{rem}
It is well known that vector $A_{1}^{\dag }b_{1}$ has the smallest Euclidean
norm among all solutions of the unconstrained least squares problem $%
\min_{x\in \mathbb{C}^{n}}\left\Vert A_{1}x-b_{1}\right\Vert $. We have shown in part 2. of Corollary~\ref{cor: new2} that $\xi = A_1^\dag b_1 + P_{\mathcal{Y},\mathcal{X}}(A_2^\dag b_2-A_1^\dag b_1)$ is the shortest solution of the constrained least squares problem (\ref{eq: cls}). 

Eld\'{e}n \cite[Theorem 2.1]{elden.80} studied minimal norm solutions of constrained least squares. On setting 
\[
h=b_2-A_2A_1^\dag b_1,\quad f=x-A_1^\dag b_1,\quad K=A_1,\quad L=A_2,\quad M=I,
\]
Eld\'{e}n's solution yields that
\[
\zeta:=A_{1}^{\dag }b_{1}+(I-P_{N(A_{2})}(A_{1}P_{N(A_{2})})^{\dag
}A_{1})A_{2}^{\dag }(b_{2}-A_{2}A_{1}^{\dag }b_{1})
\]%
minimizes the Euclidean distance $||x-A_{1}^{\dag }b_{1}||$ among all constrained minimizers $x\in \Xi $. 

With a little bit of work one finds $\zeta = \xi -P_{\mathcal{Y},\mathcal{X}}P_{N(A_2)}A_1^\dag b_1=\xi$, since 
$P_{N(A_2)}A_1^\dag\in\mathcal{X}$ by virtue of (\ref{eq: calX}). Thus part 3. of Corollary~\ref{cor: new2} simplifies and extends Eld\'{e}n's result.
\end{rem}

\end{document}